\documentclass{amsart}
\usepackage{graphicx}
\usepackage{amssymb}
\usepackage{amsmath}
\usepackage{setspace}
\newtheorem{thm}{Theorem}[]
\newtheorem{cor}[]{Corollary}
\newtheorem{lem}[]{Lemma}
\newtheorem{prop}[]{Proposition}
\theoremstyle{definition}
\newtheorem{defn}[]{Definition}

\theoremstyle{remark}

\theoremstyle{definition}

\newcommand{\tr}{{\rm tr}}

\newcommand{\Var}{{\rm Var}}
\newcommand{\la}{\lambda}

\newcommand{\calD}{\mathcal{D}}

\newcommand{\calM}{\mathcal{M}}

\newcommand{\bB}{\mathbb{B}}
\newcommand{\bT}{\mathbb{T}}

\newcommand{\bK}{\mathbb{K}}
\newcommand{\bW}{\mathbb{W}}
\newcommand{\bC}{\mathbb{C}}
\newcommand{\I}{\mathbb I}
\newcommand{\J}{\mathbb J}
\newcommand{\fA}{\mathbf{a}}

\newcommand{\mX}{\mathbf X}

\newcommand{\mA}{\mathbf A}
\newcommand{\mB}{\mathbf B}
\newcommand{\mC}{\mathbf C}
\newcommand{\mE}{\mathbf E}

\newcommand{\mI}{\mathbf I}

\newcommand{\mP}{\mathbf P}

\author{Tamer Oraby}
\thanks{Department of Mathematical Sciences,
University of Cincinnati, 2855 Campus Way, P.O. Box 210025,
Cincinnati, OH 45221-0025, USA. E-mail: orabyt@math.uc.edu. }

\title[Limiting spectral laws of random block-matrices]{The spectral laws of Hermitian block-matrices with large random blocks}

\begin{document}

\maketitle
\begin{abstract}
We are going to study the limiting spectral measure of fixed
dimensional Hermitian block-matrices with large dimensional Wigner
blocks. We are going also to identify the limiting spectral measure
when the Hermitian block-structure is Circulant. Using the limiting
spectral measure of a Hermitian Circulant block-matrix we will show
that the spectral measure of a Wigner matrix with $k-$weakly
dependent entries need not to be the semicircle law in the limit.

\end{abstract}

\section{Preliminaries and main results}\label{sec2}

Let $\calM_n(\mathbb{C})$ be the space of all $n \times n$ matrices
with complex-valued entries. Define the normalized trace of a matrix
$\mA=\left(A_{ij}\right)_{i,j=1}^{n} \in \calM_n(\mathbb{C})$ to be
$ \tr_n(\mA):=\frac1n \sum_{i=1}^n A_{ii}$.

\begin{defn}The spectral measure
of a Hermitian $n\times n$ matrix $\mA$ is the probability measure
$\mu_{\mA}$ given by
$$\mu_{\mA}=\frac{1}{n}\sum_{j=1}^n \delta_{\la_j}$$
where $\la_1\leq\la_2\leq\cdots\leq\la_n$ are the eigenvalues of
$\mA$ and $\delta_x$ is the point mass at $x$.
\end{defn}

The weak limit of the spectral measures $\mu_{\mA_n}$ of a sequence
of matrices $\{\mA_n\}$ is called the limiting spectral measure. We
will denote the weak convergence of a probability measure $\mu_n$ to
$\mu$ by $$ \mu_n \xrightarrow{\calD} \mu \mbox{ as $n\to \infty$}.
$$

\begin{defn}
A finite symmetric block-structure $\bB_k(a,b,c,\dots)$ (or shortly
$\bB_k$) over a finite alphabet $\mathcal{K}=\{a,b,c,\dots\}$ is a
$k\times k$ symmetric matrix whose entries are elements in
$\mathcal{K}$.
\end{defn}
If $\bB_k$ is a $k\times k$ symmetric block-structure and
$\mA,\mB,\mC,\dots$ are $n\times n$ Hermitian matrices, then
$\bB_k(\mA,\mB,\mC,\dots)$ is an $nk\times nk$ Hermitian matrix. One
of the interesting block structures is the $k\times k$ symmetric
Circulant over $\{a_1,a_2,\ldots,a_k\}$ that is defined as
\begin{equation}\label{circulant}
\bC_k(a_1,a_2,\ldots,a_k)= \frac{1}{\sqrt{k}}\left[
{\begin{array}{*{20}c}
   a_1 & a_2 & a_3 & {\dots} & a_k  \\
   a_k & a_1 & a_2 & {\dots} & a_{k-1}  \\
   a_{k-1} & a_k & a_1& {\dots} &a_{k-2}   \\
   {\vdots} &  {\vdots} & {\vdots} & {\ddots} &{\vdots}   \\
   a_2 & a_3& a_4 & {\dots} & a_1  \\
 \end{array} } \right]
\end{equation}
where $a_j=a_{k-j+2}$ for $j=2,3,\ldots,k$.

A random matrix $\mA$ is a matrix whose entries are random
variables. If $\bB_k$ is a block-structure and $\mA,\mB,\mC,\dots$
are random matrices, then $\bB_k(\mA,\mB,\mC,\dots)$ is a random
block-matrix.

\begin{defn}\label{wigner}
We call an $n\times n$ Hermitian random matrix
$\mA=\frac{1}{\sqrt{n}}(X_{ij})_{i,j=1}^n$ a Wigner matrix if
$\{X_{ij}; 1 \leq i < j\}$ is a family of independent and
identically distributed complex random variables such that
$E(X_{12})=0$ and $E(|X_{12}|^2)=\sigma^2$. In addition, $\{X_{ii};
i \geq 1\}$ is a family of independent and identically distributed
real random variables that is independent of the upper-diagonal
entries. We will denote all such Wigner matrices of order $n$ by
$\mathbf{Wigner}(n,\sigma^2)$.
\end{defn}

If $\{\mA_n\} $ is a sequence of $\mathbf{Wigner}(n,\sigma^2)$
matrices, then by Wigner's Theorem (\emph{cf.} \cite{Bai99}),
$$ \mu_{\mA_n}\xrightarrow{\calD} \gamma_{0,\sigma^2} \: \mbox{ as
$n\to \infty$} \qquad a.s.$$ where $\gamma_{\alpha,\sigma^2}$ is the
semicircle law centered at $\alpha$ and of variance $\sigma^2$ which
is given as
$$\gamma_{\alpha,\sigma^2}(dx)=\frac{1}{2\pi \sigma^2} \sqrt{4 \sigma^2 - (x-\alpha)^2} \:\:
\mathbf{1}_{[\alpha-2\sigma, \alpha+2\sigma]}(x) dx.$$ Now we are
ready to state the main result of this paper.

\begin{thm}[\textbf{Existence Theorem}]\label{T1}
Consider a family of independent Wigner matrices $\left (
\{\mA_n(i)\} ;i=1,\ldots,h \right )$ for which $E(|A_{12}(i)|^4)<
\infty$ and $E(A_{11}^2(i))<\infty$ for every $i$. For a fixed
$k\times k$ symmetric block-structure $\bB_k$, define $$
\mX_{n,k}:=\bB_k(\mA_n(1), \mA_n(2),\dots, \mA_n(h)). $$ Then there
exists a non-random symmetric probability measure $\mu_{\bB_k}$
which depends only the block- structure $\bB_k$ such that $$
\mu_{\mX_{n,k}}\xrightarrow{\calD} \mu_{\bB_k} \quad\mbox{as $n\to
\infty$} \qquad a.s. $$
\end{thm}

The proof of Theorem \ref{T1} relies on free probability theory and
will be given in Section \ref{Section: Proof of T1}.

Consider the symmetric Circulant block-matrix $\bC_k$ defined in
\eqref{circulant}. If $\mA_n(1),\mA_n(2),\ldots,$
\newline $\mA_n(\lfloor \frac{k}{2} \rfloor+1)$ are independent
$\mathbf{Wigner}(n,1)$ for every $n$, then Theorem \ref{T1} insures
the existence of a non-random probability measure $\nu_k$ such that
$$
\mu_{\bC_{k}(\mA_n(1),\mA_n(2),\dots,\mA_n(\lfloor \frac{k}{2}
\rfloor+1))}\xrightarrow{\calD}
 \nu_k \mbox{ as $n\to \infty$} \qquad a.s.$$
However, Theorem \ref{T1} doesn't specify $\nu_k$ but we will
identify it in the following proposition.

\begin{prop}\label{limitcirculant}
If $\mA_n(1),\mA_n(2),\ldots,\mA_n(\lfloor \frac{k}{2} \rfloor+1)$
are independent $\mathbf{Wigner}(n,1)$ for every $n$, then
$$
\mu_{\bC_{k}(\mA_n(1),\mA_n(2),\dots,\mA_n(\lfloor \frac{k}{2}
\rfloor+1))}\xrightarrow{\calD}
 \nu_k \mbox{ as $n\to \infty$} \qquad a.s.$$
where $$\nu_k=\left\{
                                       \begin{array}{ll}
                      \frac{k-1}{k}\;
\gamma_{0,\frac{k-1}{k}}+\frac{1}{k}\;\gamma_{0,\frac{2k-1}{k}}, &
\hbox{if $k$ is odd;}                   \\ \\
                     \frac{k-2}{k}\;
\gamma_{0,\frac{k-2}{k}}+\frac{2}{k}\;\gamma_{0,\frac{2k-2}{k}}, &
\hbox{if $k$ is even.}
                                       \end{array}
                                     \right.$$
\end{prop}
\begin{proof}
Since $\mA_n(j)=\mA_n(k-j+2)$ for $j=2,3,\ldots,k$; then
\cite[Theorem 3.2.2.]{davis} implies that
$\bC_{k}(\mA_n(1),\mA_n(2),\dots,\mA_n(\lfloor \frac{k}{2}
\rfloor+1))$ has the same eigenvalues as $\{\mB_n(j); \;j=1\dots,k
\}$ where
\begin{equation}\label{oddcir}
\mB_n(j):=\frac{1}{\sqrt{k}}[\mA_n(1)+2 \sum_{\ell=2}^{(k+1)/2}
\cos(\frac{2\pi (\ell-1) (j-1)}{k})\mA_n(\ell)] \end{equation} if
$k$ is odd and
\begin{equation}\label{evencir}
\mB_n(j):=\frac{1}{\sqrt{k}}[\mA_n(1)+2 \sum_{\ell=2}^{k/2}
\cos(\frac{2\pi (\ell-1) (j-1)}{k})\mA_n(\ell)+\cos((j-1)\pi)
\mA_n(\frac{k}{2}+1)] \end{equation} if $k$ is even. Hence,
$$\mu_{\bC_{k}(\mA_n(1),\mA_n(2),\dots,\mA_n(\lfloor \frac{k}{2}
\rfloor+1))}=\frac1k \sum_{j=1}^k \mu_{\mB_n(j)}.$$

Using the well known trigonometric sum $\sum_{\ell=0}^N \cos (\ell
x)=\frac12 (\frac{\sin ((N+\frac12)x)}{\sin \frac{x}{2}}+1)$, one
can easily show that
\begin{equation}\label{trigono} \sum_{\ell=0}^N \cos^2 (\ell x)=\frac12 (N+\frac32+\frac{\sin
((2N+1)x)}{\sin x}).\end{equation}

Consider the case when $k$ is odd. In Equation \eqref{oddcir}, for
$j \neq 1$, $\mB_n(j)$ is a $\mathbf{Wigner}(n,\frac{k-1}{k})$ where
the variance of the off-diagonal entries of $\mB_n(j)$ is given by
$\frac{1}{k}[1+4 \sum_{\ell=2}^{(k+1)/2} \cos^2(\frac{2\pi (\ell-1)
(j-1)}{k})]$ which turns out to be $\frac{k-1}{k}$ by Equation
\eqref{trigono}. For $j=1$, $\mB_n(1)$ is simply a
$\mathbf{Wigner}(n,\frac{2k-1}{k})$. Hence, Wigner's theorem for
$\mB_n(1)$ and the rest $k-1$ Wigner matrices $\mB_n(j);j=2,\dots,k$
finishes the proof of the odd case.

The case when $k$ is even follows from a similar argument by showing
that for $j=1,\frac{k}{2}+1$; $\mB_n(j)$ is a
$\mathbf{Wigner}(n,\frac{2k-2}{k})$ and for $j\neq 1,\frac{k}{2}+1$;
$\mB_n(j)$ is a $\mathbf{Wigner}(n,\frac{k-2}{k})$.

\end{proof}

In \cite[p.626]{Bai99}, Bai raised the question of whether Wigner's
theorem is still holding true when the independence condition in the
Wigner matrix is weakened. Schenker and Schulz-Baldes \cite{sculz}
provided an affirmative answer under some dependency assumptions in
which the number of correlated entries doesn't grow too fast and the
number of dependent rows is finite. After the first draft of the
underlying paper was completed, we learnt that Anderson and Zeitouni
\cite{AZ-2006} showed that it doesn't hold in general and they gave
an example in which the limiting spectral distribution is the free
multiplicative convolution of the semicircle law and shifted arcsine
law. In the rest of this section, we are going to use the following
corollary of Proposition \ref{limitcirculant} to give another
example.

Let $\bW(a_{11},a_{12},\ldots,a_{nn})$ be the Wigner symmetric
block-structure, \emph{i.e.},
\[\bW(a_{11},a_{12},\ldots,a_{nn})=
 \left[ {\begin{array}{*{20}c}
   a_{11} & a_{12} &  {\dots} & a_{1n}  \\
   a_{12} & a_{22} &  {\dots} & a_{2n}  \\
      {\vdots} &  {\vdots} &  {\ddots} &{\vdots}   \\
   a_{1n} & a_{2n}&  {\dots} & a_{nn}  \\
 \end{array} } \right].
\] Consider the family of $k \times k$ random matrices
$\{\mA_{ij}:i,j \geq 1 \}$ such that $\mA_{ij}=\mA_{ji}$ and
$\mA_{ij}=\bC_k(a_{ij},b_{ij},c_{ij}, \ldots)$ where
$\{a_{ij},b_{ij},c_{ij}, \ldots:i,j \geq 1  \}$ are independent and
identically distributed random variables with variance one. Then
$\bK_{n,k}:=\bW(\mA_{11},\mA_{12},\ldots,\mA_{nn})$ is an $kn \times
kn$ symmetric matrix.
\begin{cor}\label{weak} Fix $k \in \mathbb{N}$. The limiting spectral
measure of $\bK_{n,k}$ is given by
$$
\mu_{\bK_{n,k}}\xrightarrow{\calD}
 \nu_k \mbox{ as $n\to \infty$} \qquad a.s.$$
\end{cor}
In order to prove this corollary we need the following definitions.
Let $\mA$ and $\mB$ be $n\times m$ and $k\times \ell$ matrices,
respectively. By $\otimes$ we mean here the Kronecker product for
which $\mA \otimes \mB=(A_{ij}\mB)_{i=1,\ldots,n;j=1,\ldots,m}$ is
an $nk \times m \ell$ matrix. The $(p,q)$-commutation matrix
$\mP_{p,q}$ is a $pq \times pq$ matrix defined as
$$\mP_{p,q}=\sum_{i=1}^p \sum_{j=1}^q \mE_{ij} \otimes \mE_{ij}^T$$
where $\mE_{ij}$ is the $p \times q$ matrix whose entries are zero's
except the $(i,j)-$entry is 1. It is known that
$\mP_{p,q}^{-1}=\mP_{p,q}^T=\mP_{q,p}$ and $\mP_{n,k}(\mA \otimes
\mB)\mP_{\ell,m}=\mB \otimes \mA$ (\emph{cf.} \cite{searle}).
\begin{proof}
Since $\bK_{n,k}= \sum_{i,j=1}^n \widetilde{\mE}_{ij} \otimes
\mA_{ij}$ where $\widetilde{\mE}_{ij}$ is the $n \times n$ matrix
whose entries are zero's except the $(i,j)-$entry is 1. Hence
\begin{equation*}\begin{array}{l c l}
\mP_{k,n}\bK_{n,k}\mP_{n,k}&=&\sum_{i,j=1}^n \mA_{ij} \otimes
\widetilde{\mE}_{ij} \\  &=& \sum_{i,j=1}^n
\bC_k(a_{ij},b_{ij},c_{ij}, \ldots) \otimes \widetilde{\mE}_{ij} \\
&=& \sum_{i,j=1}^n
\bC_k(a_{ij}\widetilde{\mE}_{ij},b_{ij}\widetilde{\mE}_{ij},c_{ij}\widetilde{\mE}_{ij},
\ldots)  \\ &=&
 \bC_k(\sum_{i,j=1}^n a_{ij}\widetilde{\mE}_{ij}, \sum_{i,j=1}^n b_{ij}\widetilde{\mE}_{ij},\sum_{i,j=1}^n c_{ij}\widetilde{\mE}_{ij}, \ldots)  \\
&=&\bC_k(\mA_n,\mB_n,\mC_n, \ldots)
\end{array}
\end{equation*}
where $\mA_n=(a_{ij})_{i,j=1}^{ n}$, $\mB_n=(b_{ij})_{i,j=1}^{ n}$,
$\mC_n=(c_{ij})_{i,j=1}^{ n}$, $\ldots$ are independent
$\mathbf{Wigner}(n,1)$ matrices. Therefore, $\bK_{n,k}$ and
$\bC_k(\mA_n,\mB_n,\mC_n, \ldots)$ are similar to each other and so
have the same eigenvalues. Thus the result follows.
\end{proof}

Now, we define the distance on $\mathbb{N}^2$ by $d \left(
(i,j),(i',j')\right)=max \{|i-i'|,|j-j'|\} $ and for $S,T \subset
\mathbb{N}^2$; $d \left( S,T \right) =min \{d \left(
(i,j),(i',j')\right): (i,j)\in S, \, (i',j') \in T \}$. We say the
random field $\{X_{ij}:(i,j)\in \mathbb{N}^2_{\leq}\}$ is
$(k-1)$-dependent if the $\sigma$-fields
$\mathcal{F}_S=\sigma(\{X_{ij}:(i,j)\in S\})$ and
$\mathcal{F}_T=\sigma(\{X_{ij}:(i,j)\in T\})$ are independent for
all $S,T \subset \mathbb{N}^2_{\leq}$ such that $d\left( S,T
\right)>k-1$.

The matrix $\bK_{n,k}=\bW(\mA_{11},\mA_{12},\ldots,\mA_{nn})$,
defined in Corollary \ref{weak}, is an $kn \times kn$ matrix with
$(k-1)$-dependent entries, up to symmetry. That is, if we write
$\bK_{n,k}=(X_{ij})_{i,j=1}^{nk}$, then $\{X_{ij}:(i,j)\in
\mathbb{N}^2_{\leq}\}$ is a $(k-1)$-dependent random field. However,
the limiting spectral measure of $\bK_{n,k}$ is not the semicircle
law but rather a mixture of two semicircle laws due to Corollary
\ref{weak}. Our example violates the conditions imposed on the
Wigner matrix by Schenker and Schulz-Baldes in \cite{sculz} in both
the number of correlated entries and the number of dependent rows
grow as $O(n^2)$ and not $o(n^2)$.

Unfortunately, $\{X_{ij}:(i,j)\in \mathbb{N}^2_{\leq}\}$, in our
example, is not strictly stationary as the distributions remain the
same only when shifts are made by multiple of $k$.

\section{Proofs}\label{sec4}

In order to prove Theorem \ref{T1} we need to introduce some
definitions from free probability theory.

A noncommutative probability space $(\mathcal{A},\tau)$ is a pair of
a unital algebra $\mathcal{A}$ with a unit element $\I$ and a linear
functional $\tau$, called the state, for which $\tau(\I)=1$. We call
an element $\fA \in \mathcal{A}$ a noncommutative random variable
and call $\tau(\fA^n)$ its $n^{th}$ moment. We say that
$\mathcal{A}$ is a
*-algebra if the involution * is defined on $\mathcal{A}$. In
addition, we assume that $\tau(\fA^*)=\overline{\tau(\fA)}$ and
$\tau(\fA^* \fA)\geq 0$. Henceforth, we will consider only
*-algebras. We say that $\fA\in \mathcal{A} $ is selfadjoint if
$\fA^*=\fA$.

Fix a noncommutative probability space $(\mathcal{A},\tau)$. For
each selfadjoint $\fA\in \mathcal{A}$ there exists a probability
measure $\mu_\fA$ on $\mathbb{R}$ such that
$$\tau(\fA^n)=\int_{\mathbb{R}} x^n \mu_\fA(dx)$$ for all $n\geq 1$,
see \cite[p.2]{Meyer}. The probability measure $\mu_\fA$ is unique
if $|\tau(\fA^n)| \leq M^n$ for some $M>0$ and for all $n \geq 1$.

\begin{defn}[\cite{Hiai-Petz}]
A family of subalgebras $(\mathcal{A}_j; j\in J)$ of $\mathcal{A}$,
which contain $\I$, is said to be free with respect to $\tau$ if for
every $k\geq 1$ and $j_1 \neq j_2 \neq \ldots \neq j_k \in J \subset
\mathbb{N}$
$$\tau(\fA_1 \fA_2\cdots \fA_k)=0$$  for all $\fA_i \in
\mathcal{A}_{j_i}$ whenever $\tau(\fA_i)=0$ for every $1\leq i \leq
k$.
\end{defn}
Random variables in a noncommutative probability space
$(\mathcal{A},\tau)$ are said to be free if the subalgebras
generated by them and $\I$ are free.

\begin{defn}
We say that a family of sequences of random matrices
$(\{\mA_n(l)\};l=1,\ldots,m)$ is asymptotically free (\emph{cf.}
\cite{Hiai-Petz}) if for every noncommutative polynomial $p$ in $m$
variables $$\tr_n \left( p(\mA_n(1),\ldots,\mA_n(m))\right)
\xrightarrow{n\to \infty} \tau \left( p(\fA_1,\ldots,\fA_m)\right)
\qquad a.s.$$ where $(\fA_1,\ldots,\fA_m)$ is a family of free
noncommutative random variables in some noncommutative probability
space $(\mathcal{A},\tau)$.
\end{defn}

\begin{thm}\label{folk}
If $(\{\mA_n(l)\};l=1,\dots,m)$ is a family of independent
$\mathbf{Wigner}(n,1)$ matrices for which $E(|A_{12}(l)|^4)< \infty$
and $E(A_{11}^2(l))<\infty$, then $(\{\mA_n(l)\};l=1,\dots,m)$ is
asymptotically free.
\end{thm}

\subsection{Proof of Theorem \ref{folk}}

In \cite{capitaine-2005}, Capitaine and Donati-Martin showed the
asymptotic freeness for independent Wigner matrices when the
distribution of the entries is symmetric and satisfies Poincar\'{e}
inequality. Recently, Guionnet \cite{Guionnet2006} gave a proof
where she assumes that all the moments of the entries exist. Szarek
\cite{Szarek} showed us a proof for symmetric and non-symmetric
matrices with uniformly bounded entries. Szarek's proof, in brief,
is based on concentration inequalities and some tools of operator
theory. In this paper, we are going to give a combinatorial proof
for the case of Hermitian Wigner matrices with finite variance and
fourth moment of the entries.

The Schatten $p$-norm of a matrix $\mA$ is defined as
$\|\mA\|_p:=(\tr_n|\mA|^p)^{\frac1p}$ whenever $1\leq p < \infty$,
where $|\mA|=(\mA^T \mA)^{\frac12}$. The operator norm is defined as
$\|\mA\|:=\max_{1\leq i \leq n} |\la_i|$ where $\la_i$;
$i=1,2,\ldots,n$ are the eigenvalues of $\mA$. The following three
inequalities hold true;
\begin{enumerate}
\item Domination inequality \cite[p.154]{Hiai-Petz}
\begin{equation}\label{domiantion}
|\tr_n(\mA)|\leq \|\mA\|_1 \leq \|\mA\|_p \leq \|\mA\|
\end{equation}

\item H\"{o}lder's inequality \cite[p.154]{Hiai-Petz}
\begin{equation}\label{hol}
\|\mA \mB\|_r \leq \|\mA\|_p \|\mB\|_q
\end{equation}
whenever $\frac1r=\frac1p+\frac1q$ for $p,q>1$ and $r\geq 1$.

\item Generalized H\"{o}lder's inequality
\begin{equation}\label{genhol}
\|\mA{(1)} \mA{(2)} \cdots \mA{(m)}\|_1 \leq \|\mA{(1)}\|_{p_1}
\|\mA{(2)}\|_{p_2} \cdots \|\mA{(m)}\|_{p_m}
\end{equation}
where $\mA{(1)}, \mA{(2)}, \ldots, \mA{(m)}$ are $n\times n$
matrices and $\sum_{i=1}^m \frac{1}{p_i}=1$. This inequality follows
from (\ref{hol}) by induction.
\end{enumerate}

Let $\mA=\frac{1}{\sqrt{n}}(X_{i,j})_{i,j=1}^n$ be a
$\mathbf{Wigner}(n,1)$ matrix. We define
$\widetilde{\mA}=\frac{1}{\sqrt{n}}(\widetilde{X}_{i,j})_{i,j=1}^n$
to be the matrix whose off-diagonal entries are those of $\mA$
truncated by $c/\sqrt{n}$ and standardized. We will also assume that
the diagonal entries of $\widetilde{\mA}$ are zero's. In other
words,
$$\widetilde{X}_{i,j}=\left\{
    \begin{array}{ll}
    \frac{1}{\sigma(c)} \left [ X_{i,j}\mathbf{1}_{(|X_{i,j}|\leq
c)}-E(X_{i,j}\mathbf{1}_{(|X_{i,j}|\leq
c)}) \right], & \hbox{for $i<j$;} \\
      0, & \hbox{for $i=j$.}
    \end{array}
  \right.
$$
where $\mathbf{1}_{(|X_{i,j}|\leq c)}$ is equal to one if
$|X_{i,j}|\leq c$ and zero otherwise; and
$$\sigma^2(c)=E \left[ X_{i,j}\mathbf{1}_{(|X_{i,j}|\leq
c)}-E(X_{i,j}\mathbf{1}_{(|X_{i,j}|\leq c)}) \right]^2.$$ Note that
$\sigma^2(c)\to 1$ as $c\to \infty$ and
$\Var(X_{1,2}{(j)}\mathbf{1}_{(|X_{1,2}(j)|> c)}) \leq
1-\sigma^2(c)$.

The proof of Theorem \ref{folk} resembles the proof of Wigner's
theorem given in \cite{Bai99}. We will split it into a number of
lemmas.

\begin{lem}\label{prop21} If $\left (\{\mA_n(l)\};l=1,\ldots m \right)$ is a family of independent sequences of $\mathbf{Wigner}(n,1)$
matrices for which $E(|A_{12}(l)|^4)< \infty$ and
$E(A_{11}^2(l))<\infty$ for every $l$, then
\begin{equation*}\label{eq1}  \lim_{n\to \infty}
|\tr_n \left(\mA_n{(1)} \mA_n{(2)} \cdots \mA_n{(m)}\right)-\tr_n
(\widetilde{\mA}_n{(1)} \widetilde{\mA}_n{(2)} \cdots
\widetilde{\mA}_n{(m)})|=0 \qquad a.s.
\end{equation*}
\end{lem}
\begin{proof}
First, $$ \mA_n{(1)} \mA_n{(2)} \cdots
\mA_n{(m)}-\widetilde{\mA}_n{(1)} \widetilde{\mA}_n{(2)} \cdots
\widetilde{\mA}_n{(m)}=\sum_{j=1}^m \prod_{k=1}^{j-1}
\widetilde{\mA}_n{(k)} (\mA_n{(j)}-\widetilde{\mA}_n{(j)})
\prod_{l=j+1}^{m} \mA_n{(l)}$$ with the convention that
$\prod_{k=1}^{0} \widetilde{\mA}_n{(k)}=\prod_{l=m+1}^{m}
\mA_n{(l)}=\mI$. But,

\begin{equation*}
\begin{array}{l l} |\tr_n \left( \prod_{k=1}^{j-1} \widetilde{\mA}_n{(k)}
\;(\mA_n{(j)}-\widetilde{\mA}_n{(j)})\; \prod_{l=j+1}^{m} \mA_n{(l)}
\right)| &= \\ |\tr_n \left( \prod_{l=j+1}^{m} \mA_n{(l)}\;
\prod_{k=1}^{j-1} \widetilde{\mA}_n{(k)}\;
(\mA_n{(j)}-\widetilde{\mA}_n{(j)})
\right)| &\leq \\
 \| \prod_{l=j+1}^{m}  \mA_n{(l)} \; \prod_{k=1}^{j-1}
\widetilde{\mA}_n{(k)}\|_2 \; \|
\mA_n{(j)}-\widetilde{\mA}_n{(j)}\|_2
   &\leq \\ \prod_{l=j+1}^{m} \| \mA_n{(l)} \|_{2(m-1)}\; \prod_{k=1}^{j-1}
\| \widetilde{\mA}_n{(k)}\|_{2(m-1)}\; \|
\mA_n{(j)}-\widetilde{\mA}_n{(j)}\|_2
\end{array}
\end{equation*}
for all $1 \leq j \leq m$ with the convention that $\prod_{k=1}^{0}
\| \widetilde{\mA}_n{(k)}\|_p = \prod_{l=m+1}^{m} \| \mA_n{(l)}
\|_p=1$. The last two inequalities are due to the generalized
H\"{o}lder's inequality \eqref{genhol}.

It is enough to show that $$\lim_{n \to \infty} \|
\mA_n{(j)}-\widetilde{\mA}_n{(j)}\|_2 = 0 \qquad a.s. $$ for all
$j$'s, since $\lim_{n\to \infty}\| \mA_n{(l)} \|_{2(m-1)}$ and $
\lim_{n\to \infty}\| \widetilde{\mA}_n{(k)}\|_{2(m-1)} $ are finite
almost surely (\emph{cf.} \cite[Theorem 2.12]{Bai99}) for every $l$
and $k$ due to the domination inequality \eqref{domiantion} and that
$E(|A_{12}(l)|^4)< \infty$ and $E(A_{11}^2(l))<\infty$ for every
$l$. Let $\widehat{\mA}_n{(j)}:=\mA_n{(j)}-\sigma(c)
\widetilde{\mA}_n{(j)}$ or
$\widehat{X}_{r,s}{(j)}:=X_{r,s}{(j)}-\sigma(c)
\widetilde{X}_{r,s}{(j)}$ for every $r$ and $s$. Thus,
$$\|
\mA_n{(j)}-\widetilde{\mA}_n{(j)}\|_2 \leq \|
\widehat{\mA}_n{(j)}\|_2+|1-\sigma(c)| \;
\|\widetilde{\mA}_n{(j)}\|_2 $$

By definition,
$$ \| \widehat{\mA}_n{(j)}\|_2^2= \frac{1}{n^2} \sum_{r=1}^{n} \sum_{s=1}^{n}
|\widehat{X}_{r,s}{(j)}|^2=\frac{1}{n^2} \sum_{r=1}^{n}
|\widehat{X}_{r,r}{(j)}|^2+\frac{1}{n^2} \sum_{r\neq s}
|\widehat{X}_{r,s}{(j)}|^2
$$
Note that
$$\widehat{X}_{r,s}{(j)}=\left\{
  \begin{array}{ll}
    X_{r,s}{(j)}\mathbf{1}_{(|X_{r,s}(j)|> c)}
-E(X_{r,s}{(j)}\mathbf{1}_{(|X_{r,s}(j)|> c)}), & \hbox{for $r<s$;} \\
    X_{r,r}(j), & \hbox{for $r=s$.}
  \end{array}
\right. $$ Since $E(X_{1,1}^2(j))<\infty$ then $\lim_{n \to \infty}
\frac{1}{n^2} \sum_{r=1}^{n} X_{r,r}^2{(j)}=0$ almost surely due to
the Strong Law of Large Numbers (\textit{SLLN}). Once more the
\textit{SLLN} implies that
$$\lim_{n \to \infty} \frac{1}{n^2} \sum_{r\neq s}
| \widehat{X}_{r,s}{(j)}|^2
=\Var(X_{1,2}{(j)}\mathbf{1}_{(|X_{1,2}(j)|> c)}) \quad a.s.
$$
Hence, $\lim_{n\to \infty}\|
\widehat{\mA}_n{(j)}\|_2^2=\Var(X_{1,2}{(j)}\mathbf{1}_{(|X_{1,2}(j)|>
c)})$ almost surely. It is also evident that $\lim_{n\to
\infty}\|\widetilde{\mA}_n{(j)}\|_2=1$ almost surely. Therefore, for
arbitrary small $\epsilon<0$ and sufficiently large $c$,
$$\limsup_{n\to \infty} \|
\mA_n{(j)}-\widetilde{\mA}_n{(j)}\|_2 \leq
1-\sigma^2(c)+|1-\sigma(c)| <\epsilon$$ which completes the proof.
\end{proof}
Henceforth we will assume that for all $l$ the entries
$|X_{i,j}{(l)}|\leq c$ for every $i<j$ and $X_{i,i}{(l)}=0$.

\begin{lem}[\cite{Ryan98a}]\label{rate} If $\left (\{\mA_n(l)\};l=1,\ldots m \right)$ is a family of independent sequences of $\mathbf{Wigner}(n,1)$
matrices whose entries are bounded, then
\begin{equation}\label{dy} \lim_{n \to \infty} E \left(\tr_n \left(\mA_n(1)
\mA_n(2)\cdots \mA_n(m)\right)\right )=\tau \left(\fA_{1} \fA_{2}
\cdots \fA_{m}\right) \end{equation} where $\fA_{i}$'s are some free
noncommutative random variables in $(\mathcal{A},\tau)$ such that
$\fA_{i}$ has the semicircle law $\gamma_{0,1}$ for all $i$.
\end{lem}
We say that a partition $\pi=\{ B_1,\ldots,B_p \}$ of a set of
integers is non-crossing if $a<b<c<d$ is impossible for $a,c\in B_i$
and $b,d \in B_j$ when $i \neq j$. We denote the family of all
non-crossing partitions of $\{1,\ldots,k\}$ by $\text{NC}(k)$. Also
let $\text{NC}_2(k)$ be the family of all non-crossing pair
partitions which is empty unless $k$ is even. The Catalan number
$$C_k=\frac{1}{k+1}
   \left( \begin{matrix}
     2k \\
      k
    \end{matrix} \right)$$
is equal to the size of $\text{NC}(k)$ and also the size of
$\text{NC}_2(2k)$.

If $\left (\fA_l;l=1,\ldots m \right)$ is a family of free
semicircular random variables which have mean zero and variance one,
then (\emph{cf.} \cite[Equation (8)]{Ryan98a})

\begin{equation}\label{multimoment} \tau \left(\fA_{i_1} \fA_{i_2} \cdots
\fA_{i_k}\right)= \left\{
  \begin{array}{ll}
    \sum_{\pi \in \text{NC}_2(k)} \prod_{\{p,q\}\in
\pi} \mathbf{1}_{i_p=i_q}, & \hbox{if $k$ is even;} \\
    0, & \hbox{otherwise.}
  \end{array}
\right.
\end{equation}
for any $i_1,\dots,i_k \in \{1,\dots,m\}$.
\begin{lem}\label{C30} If $\left (\{\mA_n(l)\};l=1,\ldots m \right)$ is a family of independent sequences of $\mathbf{Wigner}(n,1)$
matrices whose entries are bounded and have zero diagonal entries,
then
\begin{equation*}\label{eq12}
\sum_{n=1}^\infty \Var \left( \tr_n \left(\prod_{i=1}^m
 \mA_n(i) \right)\right)<\infty
\end{equation*}
\end{lem}
\begin{proof}
It is enough to show that

\begin{equation*}
\Var \left( \tr_n\left(\prod_{i=1}^m \mA_n(i)
\right)\right)=O(n^{-2})
\end{equation*}
We will denote the number of distinct integers among
$(i_1,\ldots,i_m)$ by $\langle\langle i_1,\ldots,i_m
\rangle\rangle$.

$$\Var \left( \tr_n\left(\prod_{i=1}^m \mA_n(i) \right)\right)=E
\left( \tr_n\left(\prod_{i=1}^m \mA_n(i) \right)\right)^2-\left[ E
\left( \tr_n\left(\prod_{i=1}^m \mA_n(i) \right) \right)
\right]^2=$$

$$\frac{1}{n^{m+2}}\sum_{\I(m,n),\J(m,n)} [ E\left( \prod_{r=1}^m
X_{i_r,i_{r+1}}(r) \prod_{s=1}^m X_{j_s,j_{s+1}}(s) \right) $$
$$\qquad\qquad\qquad\qquad\qquad  -E\left( \prod_{r=1}^m
X_{i_r,i_{r+1}}(r)\right) E\left( \prod_{s=1}^m X_{j_s,j_{s+1}}(s)
\right)]$$ where $\I(m,n)=\{(i_1,\ldots,i_m):1 \leq i_1,\ldots,i_m
\leq n \}$ and $\J(m,n)=\{(j_1,\ldots,j_m):1 \leq j_1,\ldots,j_m
\leq n \}$ with the convention that $i_{m+1}=i_1$ and $j_{m+1}=j_1$.
The term under summation is zero unless:
\begin{enumerate}
\item Each one of the unordered pairs
$\left( \{i_1,i_2\},\ldots,\{i_m,i_1\} ,
\{j_1,j_2\},\ldots,\{j_m,j_1\} \right)$ appears at least twice.
\item At least one of the unordered pairs $\left( \{i_1,i_2\},\ldots,\{i_m,i_1\} \right)$ is identical to one of
the unordered pairs $\left( \{j_1,j_2\},\ldots,\{j_m,j_1\} \right)$.
\end{enumerate}

The first condition implies that $\langle\langle
i_1,\ldots,i_m,j_1,\ldots,j_m \rangle\rangle \leq m+2$. Adding the
second condition forces at least two more integers to be
replications which implies that $\langle\langle
i_1,\ldots,i_m,j_1,\ldots,j_m \rangle\rangle \leq m$. Since
$|X_{i,j}(l)|\leq c$ then

$$\Var \left( \tr_n\left(\prod_{i=1}^m \mA_n(i) \right)\right) \leq
\frac{C}{n^2}.$$
\end{proof}

\begin{proof}[Concusion of the proof of Theorem \ref{folk}]
Lemma \ref{C30} implies that the limit in \eqref{dy} is holding in
the almost sure sense due to Borel-Cantelli lemma. In other words,
if $\left (\{\widetilde{\mA}_n(l)\};l=1,\ldots m \right)$ is a
family of independent sequences of $\mathbf{Wigner}(n,1)$ matrices
whose entries are bounded, then
\begin{equation}\label{dy1} \lim_{n \to \infty} \tr_n \left(\widetilde{\mA}_n(1)
\widetilde{\mA}_n(2)\cdots \widetilde{\mA}_n(m)\right)=\tau
\left(\fA_{1} \fA_{2} \cdots \fA_{m}\right) \qquad a.s.
\end{equation}
where $\fA_{i}$'s are some free noncommutative random variables in
$(\mathcal{A},\tau)$ such that $\fA_{i}$ has the semicircle law
$\gamma_{0,1}$ for all $i$.

Now, let $\left (\{\mA_n(l)\};l=1,\ldots m \right)$ be a family of
independent sequences of $\mathbf{Wigner}(n,1)$ matrices for which
$E(|A_{12}(l)|^4)< \infty$ and $E(A_{11}^2(l))<\infty$ for every
$l$. Then by Lemma \ref{prop21} and Equation \eqref{dy1}
\begin{equation}\label{dy2} \lim_{n \to \infty} \tr_n \left(\mA_n(1)
\mA_n(2)\cdots \mA_n(m)\right)=\tau \left(\fA_{1} \fA_{2} \cdots
\fA_{m}\right) \qquad a.s. \end{equation}

Finally, since any noncommutative polynomial $p$ can be written as a
linear combination of noncommutative monomials, then $$ \lim_{n \to
\infty} \tr_n \left( p \left( \mA_n(1), \ldots, \mA_n(m) \right)
\right) =\tau \left( p \left(\fA_{1},\ldots,\fA_{m} \right) \right)
\qquad a.s.$$

\end{proof}

\subsection{Proof of Theorem \ref{T1}}\label{Section: Proof of
T1} $\newline$ \indent

Fix $k \geq 1$ and a symmetric block-structure $\bB_k$. Let $\left (
\{\mA_n(i)\} ;i=1,\ldots,h \right )$ be a family of independent
Wigner matrices such that $E(|A_{12}(i)|^4)< \infty$ and
$E(A_{11}^2(i))<\infty$ for every $i$.

Let's introduce the noncommutative probability space $(\mathcal{A}
\bigotimes \calM_k(\bC),\tau \bigotimes \tr_k)$, where $\bigotimes$
stands for the tensor product. A typical element in $\mathcal{A}
\bigotimes \calM_k(\bC)$ is a $k \times k$ matrix whose entries are
noncommutative random variables in $\mathcal{A}$. For example,
$\bB_k(\fA_1,\ldots,\fA_h)\in \mathcal{A} \bigotimes \calM_k(\bC)$
for any $\fA_1,\ldots,\fA_h \in \mathcal{A}$. The state $\tau
\bigotimes \tr_k$ is defined by $\tau \bigotimes \tr_k
(\mA)=\frac{1}{k} \sum_{i=1}^k \tau(A_{ii})$ for any $\mA \in
\mathcal{A} \bigotimes \calM_k(\bC)$.

The proof of Theorem \ref{T1} is based on the method of moments.
First, we are going to show that for every $s\in \mathbb{N}$, the
limit of $\tr_{nk}
\left(\bB_k\left(\mA_n(1),\ldots,\mA_n(h)\right)^s\right) $ exists
as $n \to \infty$, almost surely.

Fix $s \geq 1$. We can see that the trace for the $s$-power of
$\mX_{n,k}:=\bB_k \left(\mA_n(1),\ldots,\mA_n(h)\right)$ is the
trace of some noncommutative polynomial in the matrices
$\mA_n(1),\ldots,\mA_n(h)$. In other words, $$ \tr_{nk} \left(
\mX_{n,k}^s \right)=\frac1k \sum_{i=1}^k \tr_n \left( p_i\left(
\mA_n(1),\ldots,\mA_n(h) \right) \right)$$ for some noncommutative
polynomial $p_i$ and $1 \leq i \leq k$. Theorem \ref{folk} implies
that for each $i$
$$\tr_n \left( p_i\left( \mA_n(1),\ldots,\mA_n(h) \right) \right)
\to \tau \left( p_i\left( \fA_1,\ldots,\fA_h \right) \right) \mbox{
as $n\to \infty$} \qquad a.s. $$ where $\left (\fA_l;l=1,\ldots m
\right)$ is a family of free semicircular random variables.
Therefore
$$\tr_{nk} \left(\bB_k\left(\mA_n(1),\ldots,\mA_n(h)\right)^s\right)
\to \frac1k \, \tau \left( \sum_{i=1}^k  p_i\left(
\fA_1,\ldots,\fA_h \right) \right) \mbox{ as $n\to \infty$} \qquad
a.s. $$ Thus,
\begin{equation*}\label{Bconv}
\tr_{nk} \left( \mX_{n,k}^s \right) \to \tau \bigotimes \tr_k
\left(\bB_k\left(\fA_1,\ldots,\fA_h\right)^s\right) \mbox{ as $n\to
\infty$} \qquad a.s.
\end{equation*}
Note that if $s$ is an odd integer then $\tau \bigotimes \tr_k
\left(\bB_k\left(\fA_1,\ldots,\fA_h\right)^s\right)$ is zero by
Equation \eqref{multimoment}.

To complete the proof, it would be enough to show that there exist
$M> 0$ and $C >0$ such that $ \tau \bigotimes \tr_k \left(
\bB_k\left(\fA_1,\ldots,\fA_h\right)^{2s} \right) \leq C \, M^{2s}$
for all $s \geq 1$. However, for a fixed $s \geq 1$
\begin{equation*}\label{eqnbound1}
 \tau \bigotimes \tr_k \left( \bB_k\left(\fA_1,\ldots,\fA_h\right)^{2s} \right)
= \sum_{\J(2s,k)} \tau(B_{j_1 j_2} B_{j_2 j_3} \cdots B_{j_{2s}
j_1})
\end{equation*}
where $B_{ij}\in \{\fA_1,\ldots,\fA_h\}$ and
$\J(m,k):=\{(j_1,\ldots,j_m):1 \leq j_1,\ldots,j_m \leq k \}$. But
again by Equation \eqref{multimoment},
\begin{equation*}\label{eqnbound2}
  \sum_{\J(2s,k)} \tau(B_{j_1 j_2} B_{j_2 j_3} \cdots B_{j_{2s}
j_1}) \leq k^{2s}C_s=C (2k)^{2s}
\end{equation*}
for some constant $C>0$ where $C_s$ is the Catalan number.

Hence, there exists a non-random symmetric probability measure
$\mu_{\bB_k}$ with a compact support in $\mathbb{R}$ that has the
moments $\tau \bigotimes \tr_k
\left(\bB_k\left(\fA_1,\ldots,\fA_h\right)^s\right)$, for every
$s\geq 1$, such that
\begin{equation*}\label{Fconv}
\mu_{\mX_{n,k} } \xrightarrow{\calD} \,\, \mu_{\bB_k} \mbox{ as
$n\to \infty$} \qquad a.s.
\end{equation*}

\section{Concluding remarks}
\begin{enumerate}
\item We have shown in Proposition \ref{limitcirculant} that the limiting spectral measure of Hermitian
Circulant block-matrices with Wigner blocks is a mixture of two
semicircle laws. See Figure \ref{F1}.

\begin{figure}[hbt]
\begin{tabular}{ccc}
\includegraphics[height=4.5cm]{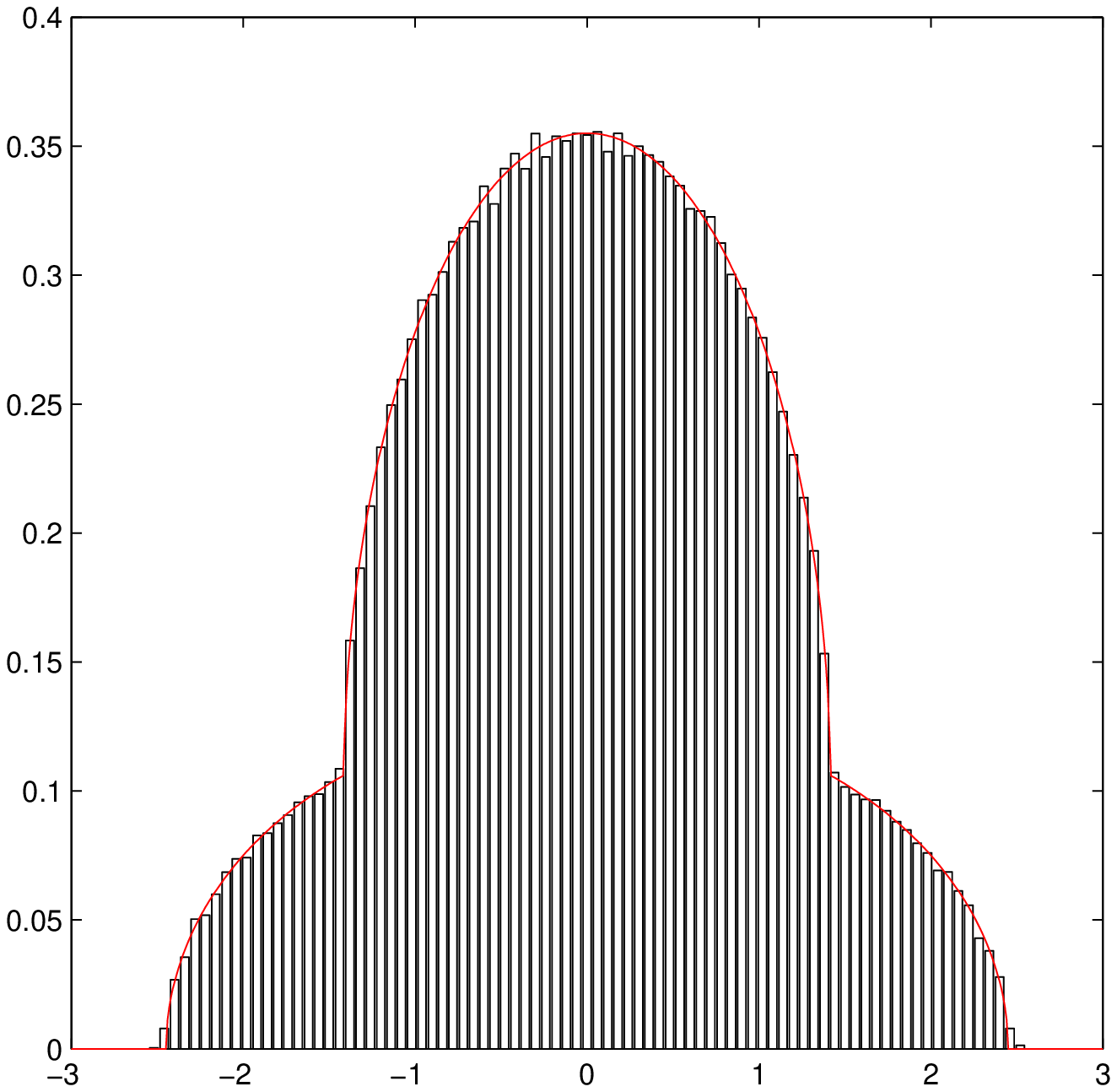}&& \includegraphics[height=4.5cm]{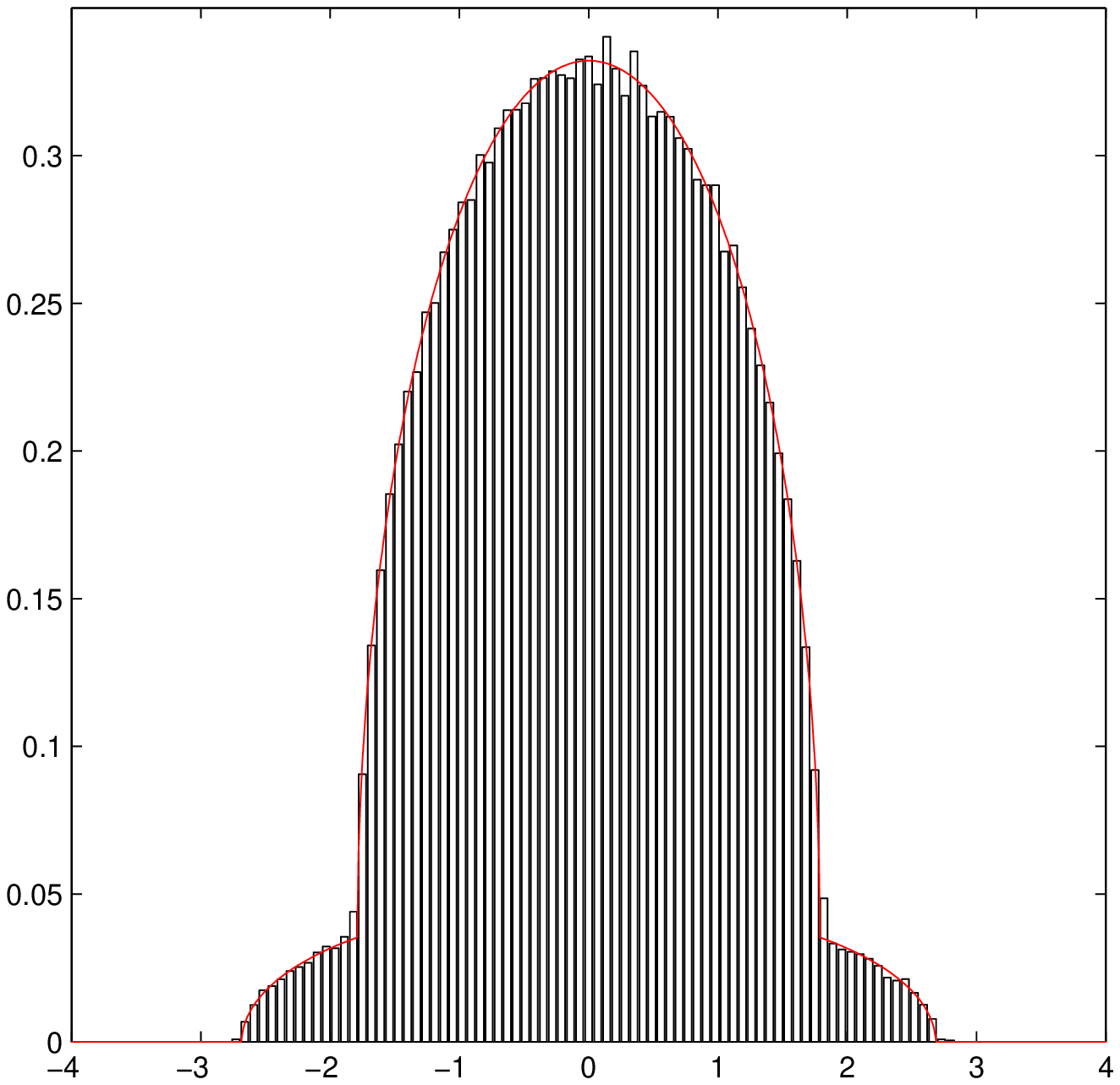}  \\
$k=4$ && $k=5$
\end{tabular}
\caption{Histograms of the eigenvalues of $100$ $\bC_k$ block
-matrices with Wigner blocks of dimension $n=200$ for $k=4$ and $5$.
The solid curves are for the exact probability density functions
provided in Proposition \ref{limitcirculant}. \label{F1}}
\end{figure}

We can also read from the simulation, see Figure \ref{F2}, of the $3
\times 3$ Toeplitz block-matrix
$$\bT_3(\mA,\mB,\mC)= \left[
{\begin{array}{*{20}c}
   \mA & \mB & \mC  \\
  \mB & \mA & \mB   \\
   \mC & \mB & \mA  \\
 \end{array} } \right]$$
that the limiting spectral measure is a mixture of two
distributions. It is evident that one of them is the semicircle law
$\gamma_{0,2}$.

\begin{figure}[hbt]
\includegraphics[height=4.5cm]{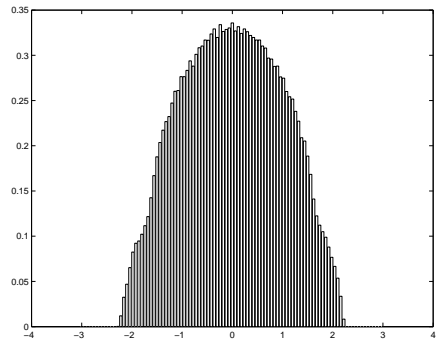}
\caption{Histograms of the eigenvalues of $100$ $\bT_3$ block
-matrices with Wigner blocks of dimension $n=200$. \label{F2}}
\end{figure}

\item If we change the blocks of the Circulant block-matrix in
Proposition \ref{limitcirculant} into random symmetric circulant
matrices then from the proof of the proposition and the limit in
\cite[Remark 2]{Bose-Mitra}, the limiting spectral measure will be a
mixture of two normal distributions.
\end{enumerate}

\section*{Acknowledgements} I would like to thank my advisor professor Bryc for his
suggestions and his continuous and sincere help and support. I am
also grateful to professor Szarek for references pertinent to
Theorem \ref{folk} and its proof.

\end{document}